\documentclass{article}
\usepackage{amsmath}
\begin{document}
	\begin{center}
		\Large{\bf\textbf\textbf{ Coupled Coincidence Point of $\phi$-Contraction Type $T$-Coupling and $(\phi,\psi)$-Contraction Type Coupling  in Metric Spaces }}
	\end{center}
	\noindent
	\begin{center}
		Tawseef Rashid \footnote{Department of Mathematics, Aligarh Muslim University, Aligarh-202002, India.\\       
			Email : tawseefrashid123@gmail.com} and Q. H. Khan \footnote{Department of Mathematics, Aligarh Muslim University, Aligarh-202002, India.\\ Email : qhkhan.ssitm@gmail.com}\\
	\end{center}
	\textbf{Abstract}. \textit{ In this research article, we discuss two topics. Firstly, we introduce SCC-Map and $\phi$-contraction type $T$-coupling. By using these two definitions, we generalize $\phi$-contraction type coupling given by H. Aydi et al. $\cite{hamba}$  to $\phi$-contraction type $T$-coupling and proved the existence theorem of coupled coincidence point for metric spaces which are not complete. Secondly, we attempt to give an answer to an open problem presented by choudhury et al. $\cite{bcm}$ concerning the investigation of fixed point and related properties for couplings satisfying other type of inequalities. In this direction we prove the existence and uniqueness theorem of strong coupled fixed point for $(\phi,\psi)$-contraction type coupling. We give  examples to illustrate our main results.}
	\noindent\\
	\textbf{2010MSC :} 47H10, 54H25, \\ 
	\begin{center}
		\textbf{1. Introduction and Mathematical Preliminaries}
	\end{center}
	T. Gnana Bhaskar and V. Lakshmikantham $\cite{gbl}$ introduced the concept of coupled fixed point of mapping $F : X \times X \to X$. Lakshmikantham V. and \'{C}iri\'{c} L. $\cite{vlc}$ introduced coupled coincidence point. Then results on existence of coupled fixed point and coupled coincidence points appeared in many papers $\cite{haezm,ha,vb,bek,bck,glk,vlc,lntx,mby,rsbm,wsbs,wsm,whnk}$ . The concept of coupling was introduced by Choudhury et al.$\cite{cm}$. Choudhury et al. $\cite{bcm}$ proved the existence and uniqueness of strong coupled fixed point for couplings using Kannan type contractions for complete metric spaces. Recently, H. Aydi et al. $\cite{hamba}$ proved the existence and uniqueness of strong coupled fixed point for $\phi$-contraction type coupling in complete partial metric spaces. In this paper we generalize result of H. Aydi et al. $\cite{hamba}$ by introducing $SCC-Map$ for metric spaces not necessarily complete. Also we prove the existence and uniqueness of coupled fixed point for $(\phi,\psi)$-contraction type coupling in complete metric spaces. Now we recall some definitions and results.
	\noindent\\  
	\textbf{Definition 1.1.} (Coupled fixed point) $\cite{gbl}$. An element $(x,y) \in X \times X$, where $X$ is any non-empty set, is called a coupled fixed point of the mapping $F : X \times X \to X$ if $F(x,y) = x~ and~ F(y,x) =y$.\\
	\textbf{Definition 1.2.} (Strong coupled fixed point) $\cite{cm}.$   An element $(x,y) \in X \times X$, where $X$ is any non-empty set,  is called a strong coupled fixed point of the mapping $F : X \times X \to X$ if $(x,y)$ is coupled fixed point and $x = y$; that is if $F(x,x) = x$.\\
	\textbf{Definition 1.3.} (Coupled Banach Contraction Mapping) $\cite{gbl}.$ Let $(X,d)$ be a metric space. A mapping $F : X \times X \to X  $ is called coupled Banach contraction if there exists $k \in (0,1)$ s.t  $\forall ~~(x,y),(u,v) \in X \times X$, the following inequality is satisfied:
	\begin{eqnarray*}
		d(F(x,y),F(u,v)) \leq \frac{k}{2} [d(x,u) + d(y,v)].\\
	\end{eqnarray*}
	\textbf{Definition 1.4.} (Cyclic mapping)  $\cite{ksv}.$ Let $A~~ and ~~B$ be two non-empty subsets of a given set $X$. Any function $f : X \to X$ is said to be cyclic (with respect to $A$ and $B$) if 
	\begin{eqnarray*}
		f(A) \subset B ~~and~~ f(B) \subset A.
	\end{eqnarray*}
	\textbf{Definition 1.5.} (Coupling) $\cite{cm}.$ Let $(X,d)$ be a metric space and $A ~and ~B$ be two non-empty subsets of $X$. Then a function $F : X \times X \to X$ is said to be a coupling with respect to $A~and~B$ if 
	\begin{eqnarray*}
		F(x,y) \in B ~and~ F(y,x) \in A~whenever~x \in A~ and~ y \in B.
	\end{eqnarray*}
	\textbf{Definition 1.6.} (Coupled coincidence point of $F$ and $g$) $\cite{vlc}.$ An element $(x,y) \in X \times X$ is called a coupled coincidence point of the mappings $F :X \times X \to X$ and $g : X \to X$ if $F(x,y) = g(x)~and~F(y,x) = g(y)$.\\
	\textbf{Definition 1.7.} \cite{ksws}.  A function $\psi: [0,\infty)\rightarrow [0,\infty)$ is called an altering distance function, if the following properties are satisfied:
	
	(i) $\psi$ is monotone increasing and continuous,
	
	(ii) $\psi(t)=0$ iff $t=0.$\\
	 \textbf{Lemma 1.8.} \cite{hamba}. Let $\phi \in \Phi$ and $\{u_n\}$ be a given sequence such that $u_n \to 0^+$ as $n \to \infty$. Then $\phi(u_n) \to 0^+$ as $n \to \infty$. Also $\phi(0) = 0$.\\
	 Where $\Phi$ is set of all functions $\phi : [0,\infty) \to [0,\infty)$ satisfying :\\
	 (i) $\phi$ is non-decreasing,\\
	 (ii) $\phi(t) < t ~~~for ~all ~~~t > 0$,\\
	 (iii) $\lim\limits_{r \to t^+}\phi(r) < t~~for ~all ~~t > 0$.\\
	 \textbf{definition 1.9.} \cite{hamba}. Let $A$ and $B$ be two non-empty subsets of a partial metric space $(X,p)$.\\
	 A coupling $F : X^2 \to X$ is said a $\phi$-contraction type coupling with respect to $A$ and $B$ if there exists $\phi \in \Phi$ such that\\
	 \begin{eqnarray*}
	 p(F(x,y),F(u,v)) \leq \phi(max\{p(x,u),p(y,v)\}),
	 \end{eqnarray*}
	 for any $x, v \in A$ and $y, u \in B$.\\
	\textbf{Theorem 1.10.} \cite{hamba}. Let $A$ and $B$ be two non-empty  closed subsets of a complete partial metric space $(X,p)$. Let $F : X^2 \to X$ be a $\phi$-contraction type coupling with respect to $A$ and $B$. Then $A \cap B \neq \emptyset$ and $F$ has a unique strong coupled fixed point in $A \cap B$.
	\begin{center}
		\textbf{2. Main Result}
	\end{center}
	\textbf{2.1. Coupled Coincidence Point of $\phi$-Contraction Type $T$-Coupling in Metric Spaces. }\\
	\noindent\\
	Before proving the main theorem of this section we introduce some definitions.\\
	\noindent\\
	\textbf{Definition 2.1.1. } (SCC-Map). Let $A$ and $B$ be any two non-empty subsets of a metric space $(X,d)$ and $T : X\to X$ be a self map on $X$. Then $T$ is said to be $SCC-map$ (with respect to $A$ and $B$), if\\
	(i) $T(A) \subseteq A$ and $T(B) \subseteq B$,\\
	(ii) $T(A)$ and $T(B)$ are closed in $X$.\\
	\noindent\\
	\textbf{Remark 2.1.2.} $I$ (the identity map) is not $SCC-Map$ in general. $I$ is $SCC-Map$ whenever $A$ and $B$ are closed subsets of $X$, i.e. $I$ can't be considered as $SCC-Map$ with respect to open sets.\\
	\noindent\\
	\textbf{Definition 2.1.3.} ($\phi$-Contraction Type $T$-Coupling). Let $A$ and $B$ be any two non-empty subsets of metric space $(X,d)$ and $T : X \to X$ is $SCC-Map$ on $X$ (w.r.t. $A$ and $B$). Then a coupling $F : X \times X \to X$ is said to be $\phi$-Contraction Type $T$-Coupling (w. r. t. $A$ and $B$), if
	\begin{eqnarray*}
	d(F(x,y),F(u,v)) \leq \phi(max\{d(Tx,Tu),d(Ty,Tv)\})
	\end{eqnarray*}
	for any $x$,$v$ $\in$ $A$ and $y$,$u$ $\in$ $B$ and $\phi \in \Phi$ defined in Lemma $1.8$.\\
	\noindent\\
	\textbf{Note 2.1.4.} If $A$ and $B$ are two non-empty subsets of a metric space $(X,d)$ and $F : X \times X \to X$ is a coupling with respect to $A$ and $B$. Then by definition of coupling for $a \in A$ and $b\in B$, we have $F(a,b) \in B$ and $F(b,a) \in A$.\\
	Now let $(a,b)$ be the coupled fixed point of $F$, then $F(a,b) = a$ and $F(b,a) = b$. But in general this is absurd because $F(a,b) \in B ~and~ a \in A$. Similarly $F(b,a) \in A ~and~ b \in B$. This is only possible for $a,b \in A \cap B$.\\
	The most important fact I want to note that for any coupling $F : X \times X \to X$ (w.r.t. $A$ and $B$), where $A$ and $B$ be any two non-empty subsets of metric space $(X,d)$, whenever we investigate for  coupled fixed point $(x,y)$ in product space  $A \times B$, then we should directly investigate in product subspace $(A \cap B) \times (A \cap B)$. Similarly for strong coupled fixed point,we should investigate  it in $A \cap B$.\\
	\noindent\\
	\textbf{Theorem 2.1.5.} Let $A$ and $B$ be any two complete subspaces of a metric space $(X,d)$ and $T : X \to X$ is $SCC-Map$ on $X$ (w.r.t. $A$ and $B$). Let $F : X \times X \to X$ be a $\phi$-contraction type $T$-coupling (with respect to $A$ and $B$), Then, 
	(i) $T(A) \cap T(B)$ $\neq$ $\emptyset$,\\
	(ii) $F$ and $T$ have atleast one coupled coincidence point in $A \times B$.\\
	\noindent\\
	\textbf{Proof.} Since $F$ is $\phi$-contractive type $T$-coupling (with respect to $A$ and $B$), we have
	\begin{equation}
		d(F(x,y),F(u,v)) \leq \phi(max\{d(Tx,Tu),d(Ty,Tv)\})
	\end{equation}
		where $x$,$v$ $\in$ $A$ and $y$,$u$ $\in$ $B$ and $\phi \in \Phi$.\\
	 As $A$ and $B$ are non-empty subsets of $X$ and $F$ is $\phi$-contraction type $T$-coupling ( w.r.t. $A$ and $B$), then for $x_0$ $\in$ $A$ and $y_0$ $\in$ $B$ we define sequences $\{x_n\}$ and $\{y_n\}$ in $A$ and $B$ resp. such that, 
	\begin{equation}
	Tx_{n+1} = F(y_n,x_n) ~~~and ~~~Ty_{n+1} = F(x_n,y_n).
	\end{equation}
	We claim $Tx_n \neq Ty_{n+1}$ and $Ty_n \neq Tx_{n+1}$   $\forall ~n$.\\
	If possible suppose for some $n$, $Tx_n = Ty_{n+1}$ and $Ty_n = Tx_{n+1}$, Then by using $(2)$, we have
	\begin{eqnarray*}
	Tx_n = Ty_{n+1} = F(x_n,y_n) ~~and ~~Ty_n = Tx_{n+1} = F(y_n,x_n).
	\end{eqnarray*}
	Which shows that $(x_n,y_n)$ is a coupled coincidence point of $F$ and $T$, so we are done in this case.
	Thus we assume 
	\begin{eqnarray*}
	 Tx_n \neq Ty_{n+1}~~ and ~~Ty_n \neq Tx_{n+1}   ~~\forall ~~n.
	\end{eqnarray*}
	Now we define a sequence $\{D_n\}$ by\\
	\begin{equation}
	D_n = max\{d(Tx_n,Ty_{n+1}),d(Ty_n,Tx_{n+1})\}.
	\end{equation}
	Then
	\begin{equation}
	D_n > 0~~~ \forall~~~ n
	\end{equation} 
	Now by using $(1)$ and $(2)$, we get
	\begin{eqnarray}
	d(Tx_n,Ty_{n+1}) &=& d(F(y_{n-1},x_{n-1}),F(x_n,y_n)\nonumber\\
	&\leq& \phi(max\{d(Ty_{n-1},Tx_n),d(Tx_{n-1},Ty_n)\})
	\end{eqnarray}
	and
	\begin{eqnarray}
	 d(Ty_n,Tx_{n+1})  &=& d(F(x_{n-1},y_{n-1}),F(y_n,x_n)\nonumber\\
	 &\leq& \phi(max\{d(Tx_{n-1},Ty_n),d(Ty_{n-1},Tx_n)\})
	\end{eqnarray}
	using $\phi(t) < t$   $\forall$   $t > 0$, then from $(3)$, $(4)$, $(5)$ and $(6)$, we have
	\begin{eqnarray}
	 0 &<& max\{d(Tx_n,Ty_{n+1}),d(Ty_n,Tx_{n+1})\}\nonumber\\
	 &\leq& \phi(max\{d(Tx_{n-1},Ty_n),d(Ty_{n-1},Tx_n)\})\\
	 &=& D_{n-1}\nonumber.
	\end{eqnarray}
	Thus $D_n < D_{n-1}$ $\forall$ $n$. This show that $\{D_n\}$ is monotonic decreasing sequence of non-negative real numbers, therefore $\exists$  $s \geq 0$, s.t. 
	\begin{equation}
	\lim\limits_{n\to+\infty}D_n = \lim\limits_{n\to+\infty}max\{d(Tx_n,Ty_{n+1}),d(Ty_n,Tx_{n+1})\} = s.
	\end{equation}
	suppose $s > 0$, letting $n \to +\infty$ in $(7)$, using $(4)$ and Lemma $1.8$, we have
	\begin{eqnarray*}
	0 < s &\leq& \lim\limits_{n \to +\infty} \phi(max\{d(Tx_{n-1},Ty_n),d(Ty_{n-1},Tx_n)\})\\
	&=& \lim\limits_{t \to s^+}\phi(t) < s.
	\end{eqnarray*}
	which is a contradiction, therefore $s = 0$.\\
	So 
	\begin{eqnarray*}
	\lim\limits_{n\to+\infty}max\{d(Tx_n,Ty_{n+1}),d(Ty_n,Tx_{n+1})\} = 0.
	\end{eqnarray*}
	i.e.
	\begin{equation}
	\lim\limits_{n\to+\infty}d(Tx_n,Ty_{n+1}) = 0~ and \lim\limits_{n\to+\infty}d(Ty_n,Tx_{n+1}) = 0.
	\end{equation}
	Now we prove $\lim\limits_{n\to\infty}d(Tx_n,Ty_n) = 0$.\\
	Let us define a sequence $\{R_n\}$ by $R_n = d(Tx_n,Ty_n)$.
	If $R_{n_0} = 0$ for some $n_0$, then $Tx_{n_0} = Ty_{n_0}$ and so $Tx_{n_0+1} = Ty_{n_0+1}$, by induction we have \\
	$d(Tx_n,Ty_n) = 0 ~\forall~ n \geq n_0.$
	Thus $\lim\limits_{n\to\infty}d(Tx_n,Ty_n) = 0.$\\
	Now we assume $R_n > 0~~\forall~~n$, then by using $(1)$,$(2)$ and definition of $\phi$, we have
	\begin{eqnarray*}
	R_n = d(Tx_n,Ty_n) &=& d(F(y_{n-1},x_{n-1}),F(x_{n-1},y_{n-1})\\
	&\leq& \phi(max\{d(Ty_{n-1},Tx_{n-1}),d(Tx_{n-1},Ty_{n-1})\})\\
	&=& \phi(d(Tx_{n-1},Ty_{n-1}))\\
	&=& \phi(R_{n-1})\\
	&<& R_{n-1}.
	\end{eqnarray*}
	Thus $\{R_n\}$ is a monotonic decreasing sequence of non-negative real numbers. Therefore $\exists ~r \geq 0$, s.t.\\
	\begin{eqnarray*}
	\lim\limits_{n\to\infty}R_n = r^+
	\end{eqnarray*}
	Assume $r > 0$ and proceeding similarly by using $\lim\limits_{z\to m^+}\phi(z) < m~~\forall~~m>0$ as above we will obtain a contradiction, so $r = 0$, i.e.\\
	\begin{equation}
	\lim\limits_{n\to\infty}d(Tx_n,Ty_n) = 0
	\end{equation}
	using triangular inequality, $(9)$ and $(10)$, we have
	\begin{equation}
	\lim\limits_{n\to\infty}d(Tx_n,Tx_{n+1}) \leq \lim\limits_{n\to\infty}[d(Tx_n,Ty_{n+1}) + d(Ty_{n+1},Tx_{n+1})] = 0.
	\end{equation}
	\begin{equation}
		\lim\limits_{n\to\infty}d(Ty_n,Ty_{n+1}) \leq \lim\limits_{n\to\infty}[d(Ty_n,Tx_{n+1}) + d(Tx_{n+1},Ty_{n+1})] = 0.
	\end{equation}
	Now we show that $\{Tx_n\}$ and $\{Ty_n\}$ are Cauchy sequences in $T(A)$ and $T(B)$ resp.
	Assume either $\{Tx_n\}$ or $\{Ty_n\}$ is not a Cauchy sequence, i.e.
	\begin{eqnarray*}
	\lim\limits_{n,m\to+\infty}d(Tx_m,Tx_n) \neq 0~~~ or~~ \lim\limits_{n,m\to+\infty}d(Ty_m,Ty_n) \neq 0.
	\end{eqnarray*}
	Then $\exists~ \varepsilon > 0$, for which we can find subsequences of integers $\{m(k)\}$ and $\{n(k)\}$ with $n(k) > m(k) > k$, s.t.
	\begin{equation}
	max\{d(Tx_{m(k)},Tx_{n(k)}),d(Ty_{m(k)},Ty_{n(k)})\} \geq \varepsilon.
	\end{equation}
	Further corresponding to $m(k)$ we can choose $n(k)$ in such a way that it is smallest integer with $n(k) > m(k)$ and satisfy $(13)$, then
	\begin{equation}
	max\{d(Tx_{m(k)},Tx_{n(k)-1}),d(Ty_{m(k)},Ty_{n(k)-1})\} < \varepsilon.
	\end{equation}
	Now by using triangular inequality, $(1)$ and $(2)$, we have
	\begin{eqnarray}
	d(Tx_{n(k)},Tx_{m(k)}) &\leq& d(Tx_{n(k)},Ty_{n(k)}) + d(Ty_{n(k)},Tx_{m(k)+1})\nonumber\\
	&& + d(Tx_{m(k)+1},Tx_{m(k)})\nonumber\\
	&=& d(Tx_{n(k)},Ty_{n(k)}) +  d(Tx_{m(k)+1},Tx_{m(k)})\nonumber\\
	&& + d(F(x_{n(k)-1},y_{n(k)-1}),F(y_{m(k)},x_{m(k)}))\nonumber\\
	&\leq&  d(Tx_{n(k)},Ty_{n(k)}) +  d(Tx_{m(k)+1},Tx_{m(k)})\nonumber\\
	&& + \phi(max\{d(Tx_{n(k)-1},Ty_{m(k)}),d(Ty_{n(k)-1},Tx_{m(k)}\})\nonumber\\
	&=& d(Tx_{n(k)},Ty_{n(k)}) +  d(Tx_{m(k)+1},Tx_{m(k)})\nonumber\\
	&& + \phi(max\{d(Tx_{m(k)},Ty_{n(k)-1}),d(Ty_{m(k)},Tx_{n(k)-1})\}).\nonumber\\
	\end{eqnarray}
	similarly by using triangular inequality, $(1)$ and $(2)$, we can show that
	\begin{eqnarray}
	d(Ty_{n(k)},Ty_{m(k)}) &\leq& d(Ty_{n(k)},Tx_{n(k)}) +  d(Ty_{m(k)+1},Ty_{m(k)})\nonumber\\
	&& + \phi(max\{d(Ty_{m(k)},Tx_{n(k)-1}),d(Tx_{m(k)},Ty_{n(k)-1})\}).\nonumber\\
	\end{eqnarray}
	from $(13)$, $(15)$ and $(16)$, we have
	\begin{eqnarray}
	\varepsilon &\leq& max\{d(Tx_{n(k)},Tx_{m(k)}),d(Ty_{n(k)},Ty_{m(k)})\}\nonumber\\ &\leq& d(Tx_{n(k)},Ty_{n(k)}) + max\{d(Tx_{m(k)+1},Tx_{m(k)}),d(Ty_{m(k)+1},Ty_{m(k)})\}\nonumber\\
	&& + \phi(max\{d(Tx_{m(k)},Ty_{n(k)-1}),d(Ty_{m(k)},Tx_{n(k)-1})\}) 
	\end{eqnarray}
	by triangular inequality , we have
	\begin{equation}
	d(Tx_{m(k)},Ty_{n(k)-1}) \leq d(Tx_{m(k)},Ty_{m(k)}) + d(Ty_{m(k)},Ty_{n(k)-1}).
	\end{equation}
	\begin{equation}
	d(Ty_{m(k)},Tx_{n(k)-1}) \leq d(Ty_{m(k)},Tx_{m(k)}) + d(Tx_{m(k)},Tx_{n(k)-1}).
	\end{equation}
	from $(14)$, $(18)$ and $(19)$, we have\\
	\noindent\\
	$max\{d(Tx_{m(k)},Ty_{n(k)-1}),d(Ty_{m(k)},Tx_{n(k)-1})\}$
	\begin{eqnarray}
	&\leq& d(Tx_{m(k)},Ty_{m(k)})
	+ max\{ d(Ty_{m(k)},Ty_{n(k)-1}),d(Tx_{m(k)},Tx_{n(k)-1})\}\nonumber\\
	& <&  d(Tx_{m(k)},Ty_{m(k)}) + \varepsilon.
	\end{eqnarray}
	since $\phi$ is non-decreasing, we have from $(20)$
	\begin{equation}
	\phi(max\{d(Tx_{m(k)},Ty_{n(k)-1}),d(Ty_{m(k)},Tx_{n(k)-1})\}) < \phi( d(Tx_{m(k)},Ty_{m(k)}) + \varepsilon).
	\end{equation}
	Now using $(21)$ in $(17)$, we get
	\begin{eqnarray}
	\varepsilon &<&  d(Tx_{n(k)},Ty_{n(k)}) + max\{d(Tx_{m(k)+1},Tx_{m(k)}),d(Ty_{m(k)+1},Ty_{m(k)})\}\nonumber\\
	&& + \phi( d(Tx_{m(k)},Ty_{m(k)}) + \varepsilon).
	\end{eqnarray}
	letting $k \to \infty$ in $(22)$ and using $(10)$, $(11)$, $(12)$ and property $(iii)$ of $\phi$ in Lemma $1.8$, we get
	\begin{eqnarray}
	\varepsilon &<& \lim\limits_{k\to\infty}\phi( d(Tx_{m(k)},Ty_{m(k)}) + \varepsilon)\nonumber\\
	&=& \lim\limits_{d(Tx_{m(k)},Ty_{m(k)}) + \varepsilon \to \varepsilon^+}\phi(d(Tx_{m(k)},Ty_{m(k)}) + \varepsilon)\nonumber\\
	&<& \varepsilon.
	\end{eqnarray}
	which is a contradiction. Thus $\{Tx_n\}$ and $\{Ty_n\}$ are Cauchy sequences in $T(A)$ and $T(B)$ resp.But $T(A)$ and $T(B)$ are closed subsets of a complete subspaces $A$ and $B$ resp. Hence  $\{Tx_n\}$ and $\{Ty_n\}$ are convergent in $T(A)$ and $T(B)$ resp.\\
	So $\exists~~ u \in T(A)$ and $v \in T(B)$, s.t.
	\begin{equation}
	Tx_n \to u ~~~and~~~ Ty_n \to v.
	\end{equation}
	 using $(10)$ in above, we get 
	 \begin{equation}
	 u = v.
	 \end{equation}
	 Therefore $u = v \in T(A) \cap T(B)$, thus $T(A) \cap T(B) \neq \emptyset$. This proves part $(i)$.
	 Now, as $u \in T(A)$ and $v \in T(B)$, therefore $\exists~a \in A ~and~ b \in B$, s.t.\\
	 $u = T(a)$ and $v = T(b)$. Using in $(24)$, we get 
	 \begin{equation}
	  Tx_n \to T(a) ~~ and ~~Ty_n \to  T(b).
	 \end{equation}
	 also from $(25)$, we get
	 \begin{equation}
	 T(a) =T(b).
	 \end{equation}
	 Now using triangular inequality, $(1)$, $(2)$, $(26)$, $(27)$ and Lemma $1.8$, we get
	 \begin{eqnarray*}
	 d(T(a),F(a,b)) &\leq& d(T(a),Tx_{n+1}) + d(Tx_{n+1},F(a,b))\\
	 &=& d(T(a),Tx_{n+1}) + d(F(y_n,x_n),F(a,b))\\
	 &\leq& d(T(a),Tx_{n+1}) + \phi(max\{d(Ty_n.T(a)),d(Tx_n,T(b))\})\\
	 &\to& 0 ~~~~as~~~~n\to \infty.
	 \end{eqnarray*}
	 thus from above, we have
	 \begin{equation}
	  F(a,b) = T(a).
	 \end{equation}
	 Again by  using triangular inequality, $(1)$, $(2)$, $(26)$, $(27)$ and Lemma $1.8$, we get
	 \begin{eqnarray*}
	 d(T(b),F(b,a)) &\leq& d(T(b),Ty_{n+1}) + d(Ty_{n+1},F(b,a))\\
	 &=& d(T(b),Ty_{n+1}) + d(F(x_n,y_n),F(b,a))\\
	 &\leq& d(T(b),Ty_{n+1}) + \phi(max\{d(Tx_n.T(b)),d(Ty_n,T(a))\})\\
	 &\to& 0 ~~~~as~~~~n\to \infty.
	 \end{eqnarray*}
	 from above, we have
	 \begin{equation}
	 F(b,a) = T(b).
	 \end{equation}
	 Hence $(28)$ and $(29)$ shows that $(a,b) \in A \times B$ is the coupled coincidence point of $F$ and $T$.\\
	 \noindent\\
	 \textbf{Corollary 2.1.6.} It should be noted that the above condition also gives a symmetric point of $F$ in $A \times B$, i.e. there exists a point $(a,b) \in A \times B$ s.t. $F(a,b) = F(b,a)$. This can be easily see by using $(27)$ in $(28)$ and $(29)$ of Theorem $2.1.5$, we get $F(a,b) = F(b,a)$.\\
	 \noindent\\
	 \textbf{Corollary 2.1.7.} If we take $T = I$ (the identity map) and $A$ and $B$ the closed subsets, then Theorem $2.1.5$ will reduce to Theorem $1.10$ by H.Aydi $\cite{hamba}$ for metric spaces not necessary complete.\\
	 \textbf{Proof :} The proof can be easily verified by using $(27)$ and $(28)$ of Theorem $2.1.5$ and the fact that $I$ is one-one map, so $a = b$ and hence $A \cap B \neq \emptyset$ and $F(a,a) =a$.\\
	 \noindent\\
	 \textbf{Note 2.1.8.} It should be noted that if $T$ is one-one, then by the asssumption $T(A) \subseteq A$ and $T(B) \subseteq B$, we have $T$ is identity map on $A$ and $B$ and uniqueness can be proved by Corollary $2.1.7.$\\
	 \noindent\\
	 \textbf{Example 2.1.9.} Let $X = (-5,5)$ be the metric space with respect usual metric $d$ on $X$, i.e. $d(x,y) = \mid x-y\mid$. Let $A = [0,2]$ and $B = [0,4]$ be the complete subspaces of $X$. Let us define $F : X \times X \to X$ by
	 \begin{equation}
	 	F(x,y) = \begin{cases}
	 		2,  & 0 \leq x, y \leq 2\\
	 		\frac{x + y}{24}, & elsewhere.\\
	 	\end{cases}
	 \end{equation}
	 Also we define $T : X \to X$ by
	 \begin{equation}
	 	T(x) = \begin{cases}
	 		2, & 0 \leq x \leq 2\\
	 		4, & x>2.\\
	 	\end{cases}
	 \end{equation}
	 We define $\phi : [0,\infty) \to [0,\infty)$ by
	 \begin{equation}
	 \phi(t) = \begin{cases}
	      \frac{2}{3}t, & 0 \leq t \leq \frac{47}{24}\\
	      \frac{47}{24}, & t>\frac{47}{24}.\\
	 \end{cases}
	 \end{equation}
	 clearly $\phi \in \Phi$.\\
	 Now from $(31)$, we have
	 \begin{eqnarray*}
	 T(A) = \{2\} \subseteq A~~ and~~ T(B) = \{2,4\} \subseteq B
	 \end{eqnarray*}
	 also $T(A)$ and $T(B)$ are closed in $A$ and $B$ resp. Thus $T : X \to X$ is SCC-Map (w.r.t. $A$ and $B$).\\
	 Now we will show $F : X \times X \to X$ is coupling (w.r.t. $A$ and $B$).\\
	 let $x \in A$ and $y \in B$.\\
	 Here two cases will arise for $y$,\\
	 case(i): $0 \leq y \leq 2$, i.e. $y \in A$,\\
	 case(ii): $2 < y \leq 4$.\\
	 For case (i) i.e. $x ,y \in A$, by using $(30)$, we have\\ $F(x,y) = 2 \in B$ and $F(y,x) = 2 \in A$. Thus $F$ is coupling (w.r.t. $A$ and $B$) and we are done in this case.\\
	 For case (ii) i.e. $ x \in A$
     and $2 < y \leq 4$, by using $(30)$, we have
     \begin{eqnarray*}
     F(x,y) = \frac{x + y}{24}
     \end{eqnarray*}
     i.e.
     \begin{eqnarray*}
     \frac{1}{12} < F(x,y) \leq \frac{1}{4}, ~~\Rightarrow F(x,y) \in B.
     \end{eqnarray*}
     and
     \begin{eqnarray*}
     	\frac{1}{12} < F(y,x) \leq \frac{1}{4}, ~~\Rightarrow F(y,x) \in A.
     \end{eqnarray*}
     thus in both the cases we get $F$ is a coupling (w.r.t. $A$ and $B$).\\
     Now we show that $F$ is $\phi$-contraction type $T$-coupling (w.r.t. $A$ and $B$.) \\
     let $x, v \in A$ and $y, u \in B$, three cases will arise for $y, u$,\\
     case(i): when both $y, u \in A$, i.e. $0 \leq y,u \leq 2.$\\
     case(ii): when one is in $A$ and other outside $A$.\\
     case(iii): when both $y,u$ lie outside $A$, i.e. $2 <y,u \leq 4.$\\
     For case(i), i.e. $x,y,u,v \in A$, we have from $(31)$\\
     
     $T(x) = T(y) = T(u) = T(v) = 2.$\\
    
    so,  $d(T(x),T(u)) = d(T(y),T(v)) = 0 $\\
    thus
    \begin{eqnarray*}
    max\{d(T(x),T(u)),d(T(y),T(v))\} = 0.
    \end{eqnarray*}
    Using $(32)$ in above, we get 
    \begin{equation}
    \phi( max\{d(T(x),T(u)),d(T(y),T(v))\}) = \phi(0) = 0.
    \end{equation}
    also for  $x,y,u,v \in A$, we have from $(30)$
    \begin{eqnarray}
    F(x,y) = F(u,v) = 2 , ~~~\Rightarrow~~~d(F(x,y),F(u,v)) = 0.
    \end{eqnarray}
    thus from $(33)$ and $(34)$, we get
    \begin{eqnarray*}
    d(F(x,y),F(u,v)) = \phi( max\{d(T(x),T(u)),d(T(y),T(v))\}).
    \end{eqnarray*}
    hence we have proved in this case.\\
    For case(ii), i.e. $x,v \in A$ and either $y$ or $u$ $\in$ A. Without loss of generality we assume $y \in A$ and $u$ outside $A$ i.e. $2 < u \leq 4$.\\
    thus for $x,y,v \in A$ and $ 2 < u \leq 4$, we have from $(31)$\\
    
   $ T(x) = T(y) = T(v) = 2 ~~and~~ T(u) = 4.$\\
   
   so ~~$d(T(x),T(u)) = 2~~ and~~d(T(y),T(v)) = 0$\\
   thus 
   \begin{eqnarray*}
   	max\{d(T(x),T(u)),d(T(y),T(v))\} = 2.
   \end{eqnarray*}
    Using $(32)$ in above, we get 
    \begin{equation}
    \phi( max\{d(T(x),T(u)),d(T(y),T(v))\}) = \phi(2) = \frac{47}{24}.
    \end{equation}
    also for  $x,y,v \in A$ and $2 < u \leq 4$, we have from $(30)$
    \begin{equation}
    F(x,y) = 2~~~ and~~~ \frac{1}{12} < F(u,v) \leq \frac{1}{4}.
    \end{equation}
    therefore from  $(36)$, we have 
    \begin{eqnarray}
    d(F(x,y),F(u,v)) &<& (2 - \frac{1}{12})\nonumber\\ &=& \frac{23}{12}\nonumber\\
    & <&  \frac{47}{24}. 
    \end{eqnarray}
    thus from $(35)$ and $(37)$, we get
    \begin{eqnarray*}
    d(F(x,y),F(u,v)) < \phi( max\{d(T(x),T(u)),d(T(y),T(v))\}).
    \end{eqnarray*}
    which proves case(ii).\\
    For case(iii), i.e. $x,v \in A$ and $2 < y,u \leq 4.$, we have from $(31)$
    \begin{eqnarray*}
    T(x) = T(v) = 2 ~and ~T(y) = T(u) = 4, 
    \end{eqnarray*}
    so, 
    \begin{eqnarray*}
     d(T(x),T(u)) = d(T(y),T(v)) = 2.
    \end{eqnarray*}
     and
     \begin{eqnarray*}
     	max\{d(T(x),T(u)),d(T(y),T(v))\} = 2.
     \end{eqnarray*}
     Using $(32)$ in above, we get 
     \begin{equation}
     \phi( max\{d(T(x),T(u)),d(T(y),T(v))\}) = \phi(2) = \frac{47}{24}.
     \end{equation}
     also for  $x,v \in A$ and $2 < y,u \leq 4$, we have from $(30)$
     \begin{equation}
      \frac{1}{12} < F(x,y) \leq \frac{1}{4} ~~~and~~~\frac{1}{12} < F(u,v) \leq \frac{1}{4}
     \end{equation}
     from $(39)$, we have
     \begin{eqnarray}
     d(F(x,y),F(u,v)) < (\frac{1}{4} - \frac{1}{12})& =& \frac{1}{6}\nonumber\\
     & < & \frac{47}{24}. 
     \end{eqnarray}
     from $(38)$ and $(40)$, we have
     \begin{eqnarray*}
     d(F(x,y),F(u,v)) < \phi( max\{d(T(x),T(u)),d(T(y),T(v))\}).
     \end{eqnarray*}
    Thus in all the cases we have proved that $F$ is $\phi$-contraction type $T$-coupling (w.r.t. $A$ and $B$).\\
    Hence all the assumptions of Theorem $2.1.5$ are satisfied, therefore $F$ and $T$ have coupled coincidence point in $A \times B$.\\
    For $a \in A$ and $b \in B ~~s.t. 0\leq b \leq 2$, then from $(30)$ and $(31)$
    \begin{equation}
    F(a,b) = 2 = T(a) ~~~and~~~F(b,a) = 2 = T(b).
    \end{equation}
    This shows that $(a,b)$ is coupled coincidence point of $F$ and $T$.\\
    The above example also shows that $F$ and $T$ have infinitely many coupled coincidence points.\\

	 \textbf{2.2 Strong Coupled Fixed Point For $(\phi,\psi)$-Contraction Type Coupling.}\\
	 \noindent\\
	 In this section we give answer to an open problem presented by choudhury et al.$\cite{bcm}$ concerning the investigation of fixed point and related properties for couplings satisfying other type of inequalities. Here we use $(\phi,\psi)$-contraction.
	 Before going to the main theorem of this section, we define $(\phi,\psi)$-contraction type coupling.\\
	 \noindent\\
	 \textbf{Definition 2.2.1.} $((\phi,\psi)$-contraction type coupling). Let $A$ and $B$ be two non-empty subsets of a metric space $(X,d)$ and $\phi$, $\psi$ are two altering distance functions. Then a coupling $F : X \times X \to X$ is said to be $(\phi,\psi)$-contraction type coupling (with respect to $A$ and $B$) if it satisfies the following inequality :
	 \begin{eqnarray*}
	 \psi (d(F(x,y),F(u,v))) \leq \psi(max\{d(x,u),d(y,v)\}) -\phi(max\{d(x,u),d(y,v)\}).
	 \end{eqnarray*}
	 for any $x,v \in A$ and $y,u \in B.$\\
	 \noindent\\
	 \textbf{Theorem 2.2.2}. Let $A$ and $B$ be two non-empty closed subsets of a complete metric space $(X,d)$ and $F : X \times X \to X$ is a $(\phi,\psi)$-contraction type coupling  ~( with respect to $A$ and $B$) i.e. there exists altering distance functions $\psi$, $\phi$ s.t.
	 	\begin{equation}
	 	\psi (d(F(x,y),F(u,v))) \leq \psi(max\{d(x,u),d(y,v)\}) -\phi(max\{d(x,u),d(y,v)\}).
	 	\end{equation}
	 	for any $x,v \in A$ and $y,u \in B,$ Then\\
	 	(i) $A \cap B \neq \emptyset$,\\
	 	(ii) $F$ has a unique strong coupled fixed point in $A \cap B$.\\
	 	\textbf{Proof :} Since $A$ and $B$ are non-empty subsets of $X$ and $F$ is a coupling ( w.r.t. $A$ and $B$), then for $x_0$ $\in$ $A$ and $y_0$ $\in$ $B$ we define sequences $\{x_n\}$ and $\{y_n\}$ in $A$ and $B$ resp. such that, 
	 	\begin{equation}
	 	x_{n+1} = F(y_n,x_n) ~~~and ~~~y_{n+1} = F(x_n,y_n).
	 	\end{equation} 
	 	If for some $n$, $x_{n+1} = y_n$ and $y_{n+1} = x_n$, then by using $(43)$, we have 
	 	\begin{eqnarray*}
	 		x_n = y_{n+1} = F(x_n,y_n) ~~and ~~y_n = x_{n+1} = F(y_n,x_n).
	 	\end{eqnarray*}
	 	This shows that $(x_n,y_n)$ is a coupled fixed point of $F$. So, we are done in this case.\\
	 	Now assume $x_{n+1} \neq y_n$ and $y_{n+1} \neq x_n$, $\forall$ $n$.\\
	 	Let us define a sequence $\{D_n\}$ by 
	 	\begin{equation}
	 	D_n = 	max\{d(x_{n+1},y_n),d(y_{n+1},x_n)\}.
	 	\end{equation}
	 	clearly $\{D_n\}$ $\subseteq$ $[0,\infty)$, $\forall~n$.\\
	 	Now using $(42)$, $(43)$ and fact that $x_n \in A$ and $y_n \in B$ $\forall~ n$, we have
	 	\begin{eqnarray}
	 	\psi(d(x_n,y_{n+1})) &=& \psi[d(F(y_{n-1},x_{n-1}),F(x_n,y_n))] \nonumber\\
	 	&=& \psi[d(F(x_n,y_n),F(y_{n-1},x_{n-1}))] \nonumber\\
	 	& \leq & \psi[max\{d(x_n,y_{n-1}),d(y_n,x_{n-1})\}]\nonumber\\ && - \phi[max\{d(x_n,y_{n-1}),d(y_n,x_{n-1})\}].
	 	\end{eqnarray}
	 	using properties of $\phi$, we have
	 	\begin{eqnarray*}
	 		\psi(d(x_n,y_{n+1})) \leq \psi[max\{d(x_n,y_{n-1}),d(y_n,x_{n-1})\}].
	 	\end{eqnarray*}
	 	Again using properties of $\psi$, we get
	 	\begin{equation}
	 	d(x_n,y_{n+1}) \leq max\{d(x_n,y_{n-1}),d(y_n,x_{n-1})\}.
	 	\end{equation}
	 	Now again using $(42)$, $(43)$ and fact that $x_n \in A$ and $y_n \in B$ $\forall~ n$, we have
	 	\begin{eqnarray}
	 	\psi(d(y_n,x_{n+1})) &=& \psi[d(F(x_{n-1},y_{n-1}),F(y_n,x_n)]\nonumber\\
	 	&\leq& \psi[max\{d(x_{n-1},y_n),d(y_{n-1},x_n)\}]\nonumber\\
	 	&& - \phi[max\{d(x_{n-1},y_n),d(y_{n-1},x_n)\}]. 
	 	\end{eqnarray}
	 	Now using properties of $\phi$ and $\psi$, we get
	 	\begin{equation}
	 	d(y_n,x_{n+1}) \leq  max\{d(x_{n-1},y_n),d(y_{n-1},x_n)\}.
	 	\end{equation}
	 	by using $(46)$ and $(48)$, we get
	 	\begin{eqnarray*}
	 		max\{d(y_n,x_{n+1}),d(x_n,y_{n+1})\} \leq max\{d(x_n,y_{n-1}),d(y_n,x_{n-1})\}.
	 	\end{eqnarray*}
	 	i.e.
	 	\begin{equation}
	 	max\{d(x_{n+1},y_n),d(y_{n+1}x_n)\} \leq max\{d(x_n,y_{n-1}),d(y_n,x_{n-1})\}.
	 	\end{equation}
	 	from $(44)$, we have 
	 	\begin{eqnarray*}
	 		D_n \leq D_{n-1} ~~\forall~~ n \geq 1.
	 	\end{eqnarray*}
	 	Therefore $\{D_n\}$ is monotonic decreasing sequence of non-negative real numbers.\\
	 	Thus $\exists$ $r \geq 0$, s.t. $\lim\limits_{n \to \infty}D_n = r$,\\
	 	i.e.
	 	\begin{equation}
	 	\lim\limits_{n \to \infty}max\{d(x_{n+1},y_n),d(y_{n+1},x_n)\} = r.
	 	\end{equation}
	 	Since $\psi : [0,\infty) \to [0,\infty)$ is non-decreasing, then $\forall~ a,b\in [0,\infty)$, we have
	 	\begin{equation}
	 	max\{\psi(a),\psi(b)\} = \psi(max\{a,b\}).
	 	\end{equation}
	 	on using $(45)$, $(47)$ and $(51)$, we get
	 	\begin{eqnarray*}
	    \psi[max\{d(x_n,y_{n+1}),d(y_n,x_{n+1})\}]&=& max\{\psi(d(x_n,y_{n+1}),\psi(d(y_n,x_{n+1}))\}\\
	 		&\leq& \psi[max\{d(x_n,y_{n-1}),d(y_n,x_{n-1})\}]\\
	 		&& - \phi[max\{d(x_n,y_{n-1}),d(y_n,x_{n-1})\}].
	 	\end{eqnarray*}
	 	letting $n \to \infty$ in above inequality, using$(50)$ and continuities of $\phi$ and $\psi$, we have
	 	\begin{eqnarray*}
	 		\psi(r) \leq \psi(r) - \phi(r) \leq \psi(r)
	 	\end{eqnarray*}
	 	$\Rightarrow \phi(r) = 0$, since $\phi$ is altering distance function, so $r = 0$.\\
	 	Hence $\lim\limits_{n \to \infty}D_n = 0$, i.e.
	 	\begin{eqnarray*}
	 		\lim\limits_{n \to \infty}max\{d(x_n,y_{n+1}),d(y_n,x_{n+1})\} = 0.
	 	\end{eqnarray*}
	 	Thus both
	 	\begin{equation}
	 	\lim\limits_{n \to \infty}d(x_n,y_{n+1}) = 0~~ and~~ \lim\limits_{n \to \infty}d(y_n,x_{n+1}) = 0
	 	\end{equation}
	 	Now we define a sequence $\{R_n\}$ by $R_n = d(x_n,y_n)$\\
	 	we show that $R_n \to 0~~as~~n \to \infty $.\\
	 	By using $(42)$ and $(43)$, we get
	 	\begin{eqnarray}
	 	\psi(R_n) &=& \psi(d(x_n,y_n))\nonumber\\
	 	&=&\psi(d(F(y_{n-1},x_{n-1}),F(x_{n-1},y_{n-1}))\nonumber\\
	 	&\leq& \psi[max\{d(x_{n-1},y_{n-1}),d(x_{n-1},y_{n-1})\}]\nonumber\\
	 	&&- \phi[max\{d(x_{n-1},y_{n-1}),d(x_{n-1},y_{n-1})\}]\nonumber\\
	 	&=& \psi(d(x_{n-1},y_{n-1})) - \phi((d(x_{n-1},y_{n-1}))).
	 	\end{eqnarray}
	 	by properties of $\phi$ and $\psi$, we have
	 	\begin{eqnarray*}
	 		R_n \leq d(x_{n-1},y_{n-1}) = R_{n-1},
	 	\end{eqnarray*}
	 	i.e. \begin{eqnarray*}
	 		R_n \leq R_{n-1} ~~\forall~~ n \geq 1.
	 	\end{eqnarray*}
	 	Thus $\{R_n\}$ is monotonic decreasing sequence of non-negative real numbers.\\
	 	therefore $\exists~ s \geq 0$, s.t. 
	 	\begin{eqnarray}
	 	\lim\limits_{n\to\infty}R_n = \lim\limits_{n\to\infty}d(x_n,y_n) = s
	 	\end{eqnarray}
	 	take $n \to \infty$ in $(53)$, using $(54)$ and continuities of $\phi$ and $\psi$, we have
	 	\begin{eqnarray*}
	 		\psi(s) \leq \psi(s) - \phi(s) \leq \psi(s)
	 	\end{eqnarray*}
	 	$\Rightarrow$ $\phi(s) = 0$, but since $\phi$ is altering distance function, so we have $s = 0$.\\
	 	i.e.
	 	\begin{equation}
	 	\lim\limits_{n\to\infty}R_n = \lim\limits_{n\to\infty}d(x_n,y_n) = 0.
	 	\end{equation}
	 	Now using triangular inequality, $(52)$ and $(55)$, we have
	 	\begin{equation}
	 	\lim\limits_{n\to\infty}d(x_n,x_{n+1}) \leq  \lim\limits_{n\to\infty}[d(x_n,y_n) + d(y_n,x_{n+1})] = 0
	 	\end{equation}
	 	and
	 	\begin{equation}
	 	\lim\limits_{n\to\infty}d(y_n,y_{n+1}) \leq  \lim\limits_{n\to\infty}[d(y_n,x_n) + d(x_n,y_{n+1})] = 0.
	 	\end{equation}
	 	Now we will prove that sequences $\{x_n\}$ and $\{y_n\}$ are cauchy sequences in $A$ and $B$ resp.\\
	 	If possible, let $\{x_n\}$ or $\{y_n\}$ is not a Cauchy sequence. Then there exists an $\varepsilon > 0$, and sequence of positive integers $\{m(k)\}$ and $\{n(k)\}$ such that $\forall$ positive integer $k$, with $n(k) > m(k) > k$, we have
	 	\begin{equation}
	 	T_k = max\{d(x_{m(k)},x_{n(k)}),d(y_{m(k)},y_{n(k)})\} \geq \varepsilon.
	 	\end{equation}
	 	and
	 	\begin{equation}
	 	max\{d(x_{m(k)},x_{n(k)-1}),d(y_{m(k)},y_{n(k)-1})\} < \varepsilon.
	 	\end{equation} 
	 	Now we show that 
	 	\begin{eqnarray*}
	 		d(y_{n(k)},x_{m(k)+1}) \leq 	max\{d(x_{m(k)},y_{n(k)-1}),d(y_{m(k)},x_{n(k)-1})\}
	 	\end{eqnarray*} 
	 	By using $(42)$ and $(43)$, we get
	 	\begin{eqnarray*}
	 		\psi[d(y_{n(k)},x_{m(k)+1})] &=& \psi[d(F(x_{n(k)-1},y_{n(k)-1}),F(y_{m(k)},x_{m(k)}))]\\
	 		&\leq& \psi[max\{d(x_{n(k)-1},y_{m(k)}),d(y_{n(k)-1},x_{m(k)})\}] \\
	 		&& - \phi [max\{d(x_{n(k)-1},y_{m(k)}),d(y_{n(k)-1},x_{m(k)})\}].
	 	\end{eqnarray*}
	 	Using properties of $\phi$ and $\psi$, we have
	 	\begin{equation}
	 	d(y_{n(k)},x_{m(k)+1}) \leq max\{d(x_{n(k)-1},y_{m(k)}),d(y_{n(k)-1},x_{m(k)})\}.
	 	\end{equation}
	 	Similarly we can show by same pattern that,
	 	\begin{equation}
	 	d(x_{n(k)},y_{m(k)+1}) \leq max\{d(y_{n(k)-1},x_{m(k)}),d(x_{n(k)-1},y_{m(k)})\}.
	 	\end{equation}
	 	From $(60)$ and $(61)$, we have
	 	\begin{eqnarray}
	 	max\{d(y_{n(k)},x_{m(k)+1}),d(x_{n(k)},y_{m(k)+1})\} &\leq& max\{d(x_{m(k)},y_{n(k)-1}),d(y_{m(k)},x_{n(k)-1})\}\nonumber\\
	 	& =& \lambda.
	 	\end{eqnarray}
	 	Where $\lambda = max\{d(x_{m(k)},y_{n(k)-1}),d(y_{m(k)},x_{n(k)-1})\}. $\\
	 	\noindent\\
	 	It is fact that for $a,b,c \in R^+$, $max\{a+c,b+c\} = c + max\{a,b\}$.\\
	 	Therefore by triangular inequality, $(59)$ and the above fact, we have
	 	\begin{eqnarray}
	 	\lambda &=& max\{d(x_{m(k)},y_{n(k)-1}),d(y_{m(k)},x_{n(k)-1})\}\nonumber\\
	 	& \leq&  max\{d(x_{m(k)},x_{n(k)-1}) + d(x_{n(k)-1},y_{n(k)-1}),d(y_{m(k)},y_{n(k)-1}) + d(y_{n(k)-1},x_{n(k)-1})\}\nonumber\\
	 	&=&  d(x_{n(k)-1},y_{n(k)-1}) + max\{d(x_{m(k)},x_{n(k)-1}),d(y_{m(k)},y_{n(k)-1})\}\nonumber\\
	 	&<&  d(x_{n(k)-1},y_{n(k)-1}) + \varepsilon.
	 	\end{eqnarray}
	 	Thus from $(62)$ and $(63)$, we get
	 	\begin{equation}
	 	max\{d(y_{n(k)},x_{m(k)+1}),d(x_{n(k)},y_{m(k)+1})\} <  d(x_{n(k)-1},y_{n(k)-1}) + \varepsilon.
	 	\end{equation}
	 	Now again by traingular inequality, we have
	 	\begin{equation}
	 	d(x_{n(k)},x_{m(k)}) \leq d(x_{n(k)},y_{n(k)}) + d(y_{n(k)},x_{m(k)+1}) + d(x_{m(k)+1},x_{m(k)}).
	 	\end{equation}
	 	and
	 	\begin{equation}
	 	d(y_{n(k)},y_{m(k)}) \leq d(y_{n(k)},x_{n(k)}) + d(x_{n(k)},y_{m(k)+1}) + d(y_{m(k)+1},y_{m(k)}).
	 	\end{equation} 
	 	From $(58)$ $(64)$, $(65)$ and $(66)$, we get
	 	\begin{eqnarray}
	 	T_k &=& max\{	d(x_{n(k)},x_{m(k)}),d(y_{n(k)},y_{m(k)})\}\nonumber\\
	 	&\leq& d(x_{n(k)},y_{n(k)}) + max\{d(x_{m(k)},x_{m(k)+1}),d((y_{m(k)},y_{m(k)+1})\}\nonumber\\
	 	&& + max\{d(y_{n(k)},x_{m(k)+1}),d(x_{n(k)},y_{m(k)+1)}\}\nonumber\\
	 	&<& d(x_{n(k)},y_{n(k)}) + max\{d(x_{m(k)},x_{m(k)+1}),d((y_{m(k)},y_{m(k)+1})\}\nonumber\\
	 	&& + d(x_{n(k)-1},y_{n(k)-1}) + \varepsilon. 
	 	\end{eqnarray}
	 	Take $k \to \infty$ in $(67)$ and using $(55)$, $(56)$, $(57)$ and $(58)$, we have
	 	\begin{eqnarray*}
	 		\varepsilon \leq T_k < \varepsilon.
	 	\end{eqnarray*}
	 	Which is a contradiction, Hence $\{x_n\}$ and $\{y_n\}$ are Cauchy sequences in $A$ and $B$ respectively. \\
	 	Since $A$ and $B$ are closed subsets of complete metric space X, therefore $\{x_n\}$ and $\{y_n\}$ are convergent in $A$ and $B$ respec.\\
	 	Thus $\exists~ x \in A$ and $y \in B$, s.t.
	 	\begin{equation}
	 	x_n \to x~~ and ~~y_n \to y.
	 	\end{equation}
	 	Also by using $(55)$ in $(68)$, we get
	 	\begin{equation}
	 	x = y.
	 	\end{equation}
	 	Thus $x = y \in A \cap B$, which shows that $A \cap B \neq \emptyset$.\\
	 	using $(42)$, $(43)$ and fact that $x_n,x \in A$ and $y_n,y \in B$ $\forall~ n $, we have
	 	\begin{eqnarray*}
	 		\psi(d(x_{n+1},F(x,y))) &=& \psi(d(F(y_n,x_n),F(x,y))\\
	 		&=& \psi(d(F(x,y),F(y_n,x_n))\\
	 		&\leq& \psi[max\{d(x,y_n),d(y,x_n)\}]- \phi [max\{d(x,y_n),d(y,x_n)\}].	
	 	\end{eqnarray*}
	 	using properties of $\psi$ and $\phi$, we get
	 	\begin{equation}
	 	d(x_{n+1},F(x,y)) \leq max\{d(x,y_n),d(y,x_n)\}.
	 	\end{equation}
	 	Now using triangular inequality, $(68)$, $(69)$ and $(70)$, we get
	 	\begin{eqnarray*}
	 		d(x,F(x,x)) &=& d(x,F(x,y))\\
	 		&\leq& d(x,x_{n+1}) + d(x_{n+1},F(x,y))\\
	 		&\leq& d(x,x_{n+1}) + max\{d(x,y_n),d(y,x_n)\}\\
	 		&\to& 0 ~~~as ~~~n \to \infty.
	 	\end{eqnarray*}
	 	$\Rightarrow$ $F(x,x) = x$, hence $F$ has a strong coupled fixed point in $A \cap B$.\\
	 	\textbf{Uniqueness}: Let if possible $F$ has two strong fixed points $l$, $m$ in $A \cap B$, then
	 	\begin{equation}
	 	F(l,l) = l ~~~and~~~F(m,m) = m.
	 	\end{equation}
	 	Now using $(42)$, $(71)$, properties of $\phi$ and $\psi$ and $l, m \in A \cap B$, we have
	 	\begin{eqnarray*}
	 		\psi(d(l,m)) &=& \psi(d(F(l,l),F(m,m))\\
	 		&\leq& \psi[max\{d(l,m),d(l,m)\}] - \phi [max\{d(l,m),d(l,m)\}]\\
	 		&\leq& \psi (d(l,m)) - \phi(d(l,m)).
	 	\end{eqnarray*}
	 	which gives $\phi(d(l,m)) = 0$, as $\phi$ is altering distance function so $d(l,m) = 0$.\\ Hence $l = m$, which proves the uniqueness.\\
	 	\noindent\\
	 	\textbf{Example 2.2.3. :} Let $X = [0,3]$ be the complete metric space with respect to usual metric 'd' on $X$ i.e. $d(x,y) = \mid x-y\mid$. Let $A = \{1\}$ and $B = \{1,2\}$ be the closed subsets of $X$. We define $F : X \times X \to X$ by\\
	 	\begin{equation}
	 	F(x,y) = min\{x,y\}, \forall x,y \in X.
	 	\end{equation}
	 	Also we define $\phi, \psi : [0,\infty) \to [0,\infty)$ by
	 	\begin{equation}
	 	\phi(t) = t^2 ~~~and ~~~\psi(t) = t.
	 	\end{equation}
	 	Then clearly $\phi~~and~~ \psi$ are altering distance functions and $A$ and $B$ are closed subsets of a complete metric space $[0,3].$\\
	 	First we show that $F$ is a coupling (w.r.t. $A$ and $B$).\\
	    Let $x \in A$, $y \in B$, i.e. $x =1~~ and~~ y = 1,2$, we have by $(72)$
	    \begin{eqnarray*}
	    F(x,y) = 1 \in B ~~~and~~~ F(y,x) = 1 \in A.
	    \end{eqnarray*}
	    This shows that $F$ is a coupling (w.r.t. $A$ and $B$).\\
	    Now we show that $F$ is $(\phi,\psi)$-contraction type coupling.\\
	 	let $x, v \in A$ and $y, u \in B$, then four cases arise\\
	 	case(i) $x = v = 1$ and $y = u = 1,$\\
	 	case(ii) $x = v = 1$ and $y = 1, u = 2,$\\
	 	case(iii) $x = v = 1$ and $y = 2, u = 1,$\\
	 	case(iv)  $x = v = 1$ and $y = 2, u = 2.$\\
	 	For case(i) when  $x = v = 1$ and $y = u = 1$, we have from $(72)$\\
	 	$F(x,y) = F(1,1) =1 ~~and~~ F(u,v) = F(1,1) = 1.$\\
	 	using above and $(73)$, we get
	 	\begin{equation}
	 	\psi(d(F(x,y),F(u,v))) = \psi(0) = 0.
	 	\end{equation}
	 	also $d(x,u) = 0 ~~~and~~~d(y,v) = 0,$ so $max\{d(x,u),d(y,v)\} = 0.$\\
	 	then from $(73)$, we get
	 	\begin{equation}
	 	\psi(max\{d(x,u),d(y,v)\}) - \phi(max\{d(x,u),d(y,v)\}) = \psi(0) - \phi(0) = 0.
	 	\end{equation}
	 	thus $(74)$ and $(75)$ gives
	 	\begin{eqnarray*}
	 	\psi(d(F(x,y),F(u,v))) \leq 	\psi(max\{d(x,u),d(y,v)\}) - \phi(max\{d(x,u),d(y,v)\}).
	 	\end{eqnarray*}
	 	we are done in case(i).\\
	 	For case(ii) when  $x = v = 1$ and $y =1, u = 2$, we have from $(72)$\\
	 	$F(x,y) = F(1,1) =1 ~~and~~ F(u,v) = F(2,1) = 1.$\\
	 	using above and $(73)$, we get
	 	\begin{equation}
	 	\psi(d(F(x,y),F(u,v))) = \psi(0) = 0.
	 	\end{equation}
	 	also $d(x,u) = 1 ~~~and~~~d(y,v) = 0,$ so $max\{d(x,u),d(y,v)\} = 1.$\\
	 	then from $(73)$, we get
	 	\begin{equation}
	 	\psi(max\{d(x,u),d(y,v)\}) - \phi(max\{d(x,u),d(y,v)\}) = \psi(1) - \phi(1) = 0.
	 	\end{equation}
	 	thus $(76)$ and $(77)$ gives
	 	\begin{eqnarray*}
	 		\psi(d(F(x,y),F(u,v))) \leq 	\psi(max\{d(x,u),d(y,v)\}) - \phi(max\{d(x,u),d(y,v)\}).
	 	\end{eqnarray*}
	 	which proves in case(ii).\\
	 	For case(iii) when  $x = v = 1$ and $y =2, u = 1$, we have from $(72)$\\
	 		$F(x,y) = F(1,2) =1 ~~and~~ F(u,v) = F(1,1) = 1.$\\
	 		using above and $(73)$, we get
	 		\begin{equation}
	 		\psi(d(F(x,y),F(u,v))) = \psi(0) = 0.
	 		\end{equation}
	 		also $d(x,u) = 0 ~~~and~~~d(y,v) = 1,$ so $max\{d(x,u),d(y,v)\} = 1.$\\
	 		then from $(73)$, we get
	 		\begin{equation}
	 		\psi(max\{d(x,u),d(y,v)\}) - \phi(max\{d(x,u),d(y,v)\}) = \psi(1) - \phi(1) = 0.
	 		\end{equation}
	 		thus $(78)$ and $(79)$ gives
	 		\begin{eqnarray*}
	 			\psi(d(F(x,y),F(u,v))) \leq 	\psi(max\{d(x,u),d(y,v)\}) - \phi(max\{d(x,u),d(y,v)\}).
	 		\end{eqnarray*}
	 		case(iii) is proved.\\
	 		For case(iv) when  $x = v = 1$ and $y =2, u = 2$, we have from $(72)$\\
	 		$F(x,y) = F(1,2) =1 ~~and~~ F(u,v) = F(2,1) = 1.$\\
	 		using above and $(73)$, we get
	 		\begin{equation}
	 		\psi(d(F(x,y),F(u,v))) = \psi(0) = 0.
	 		\end{equation}
	 		also $d(x,u) = 1 ~~~and~~~d(y,v) = 1,$ so $max\{d(x,u),d(y,v)\} = 1.$\\
	 		then from $(73)$, we get
	 		\begin{equation}
	 		\psi(max\{d(x,u),d(y,v)\}) - \phi(max\{d(x,u),d(y,v)\}) = \psi(1) - \phi(1) = 0.
	 		\end{equation}
	 		thus $(80)$ and $(81)$ gives
	 		\begin{eqnarray*}
	 			\psi(d(F(x,y),F(u,v))) \leq 	\psi(max\{d(x,u),d(y,v)\}) - \phi(max\{d(x,u),d(y,v)\}).
	 		\end{eqnarray*}
	 		Hence $F$ is $(\phi,\psi)$-contraction type coupling ( w.r.t $A$ and $B$). Thus all the conditions of Theorem $2.2.2.$ are satisfied:\\
	 		then $F$ has a strong coupled fixed point in $A \cap B$.\\
	 		Clearly $A \cap B = \{1\} \neq \emptyset$ and $1$ is the unique strong coupled fixed point of $F$ in $A \cap B$  as $F(1,1) = min\{1,1\} = 1$.\\


\begin{thebibliography}{23}
	\bibitem{haezm} H. Aydi, E. Karapinar, Z. Mustafa, \textit{Coupled coincidence point results on generalized distance in ordered cone metric spaces}, Positivity, \textbf{17}(4)(2013), 979-993.
	\bibitem{ha} H. Aydi, \textit{Coupled fixed point results in ordered partial metric spaces}, Sel\c{c}uk J. Appl. Math. \textbf{13}(2012),  23-33.
	\bibitem{hamba} H. Aydi, M. Barakat, A. Felhi, H. Isik, \textit{On $\phi$-contraction type couplings in partial metric spaces}, Journal Of Mathematical Analysis, \textbf{8}(4)(2017), 78-89.
	\bibitem{vb} V. Berinde, \textit{Generalized coupled fixed point theorems in partially ordered metric spaces and applications}, Nonlin. Anal. \textbf{74}(2011), 7347-7355.
	\bibitem{bek} N. Bilgili, I. M. Erham, E. Karapinar, D. Turkoglu, \textit{A note on coupled fixed point theorems for mixed g-monotone mappings in partially ordered metric spaces}, Fixed Point Theory Appl. \textbf{2014}(2014):120, 6 pages.
	\bibitem{cm} B.S. Choudhury, P. Maity,  \textit{Cyclic coupled fixed point result using Kannan type contractions}, Journal Of Operators. \textbf{2014}(2014), Article ID 876749, 5 pages.
	\bibitem{bcm}  B.S. Choudhury, P. Maity, P. Konar, \textit{Fixed point results for couplings on metric spaces}, U. P. B. Sci. Bull., series A, \textbf{79}(1)(2017). 1-12.
	\bibitem{bck} B.S. Choudhury, A. Kundu, \textit{A coupled coincidence result in partially ordered metric spaces for compatible mappings}, Nonlinear Anal. \textbf{73}(2010), 2524-2531.
	\bibitem{gbl} T. Gnana Bhaskar, V. Lakshmikantham, \textit{Fixed point theorems in partially ordered metric spaces and applications}, Nonlin. Anal. \textbf{65}(2006), 1379-1393.
	\bibitem{glk} D. Guo, V. Lakshmikantham, \textit{Coupled fixed points of non-linear operators with applications}, Nonlin. Anal. \textbf{11}(1987), 623-632.
	\bibitem{ksv} W.A. Kirk, P.S. Srinivasan and P. Veeramani, \textit{Fixed points for mappings satisfying cyclical contractive conditions},  Fixed Point Theory. \textbf{4}(2003), 78-89.
	\bibitem{ksws} M.S. Khan, M. Swaleh, S. Sessa, \textit{Fixed point theorems by altering distances between the points}, Bull. Aust. Math. Soc. \textbf{30}(1984), pp. 1-9, 1984.
	\bibitem{vlc} V. Lakshmikantham, L. \'{C}iri\'{c},  \textit{Coupled fixed point theorems for nonlinear contractions in partially ordered metric spaces}, Nonlinear Anal. \textbf{70}(2009), 4341-4349.
	\bibitem{lntx} NV. Luong, X. Thuan, \textit{Coupled fixed point theorems in partially ordered metric spacesdepended on another function}, Bull. Math. Anal. \textbf{3}(2011), 129-140.
	\bibitem{mby} Mujahid Abbas, Bashir Ali, Yusuf I.Suleiman, \textit{Generalized coupled common fixed point results in partially ordered A-metric spaces}, Fixed Point Theory and Application, (2015) 2015:64, 24 pages.
	\bibitem{rsbm} SH. Rosouli, M. Bahrampour, \textit{A remark on the coupled fixed point theorems for mixed monotone operators in partially ordered metric spaces}, J. Math. Comput. Sci. \textbf{3}(2011), 246-261.
	\bibitem{wsbs} W. Shatanawi, B. Samet, M. Abbas, \textit{Coupled fixed point theorems for mixed monotone mappings in ordered partial metric spaces}, Mathematical and Computer Modelling, \textbf{55}(2012), 680-687.
	\bibitem{wsm} Wasfi Shatanawi, Mujahid Abbas, Hassen Aydi, Nedal Tahat, \textit{Common coupled coincidence and coupled fixed points in G-metric spaces}, Nonlinear Anal. and Application, \textbf{2012}(2012), Article ID jnaa-00162, 16 pages.
	\bibitem{whnk} W. Shatanawi, H.K. Nashine, N. Tahat, \textit{Generalization of some coupled fixed point results on partial metric spaces}, International Journal of Mathematics and Mathematical Sciences, \textbf{2012}(2012), Article ID 686801, 10 pages.
	

	


	
		
	\end{thebibliography}
\end{document}